\documentclass[a4paper ]{article}
\usepackage{amsthm}
\usepackage{amsmath,amssymb}
\usepackage{amsfonts}
\usepackage{graphicx}
\usepackage{url}

\newtheorem{theorem}{Theorem}[section]
\theoremstyle{plain}
\newtheorem{conjecture}[theorem]{Conjecture}

\newtheorem{problem}[theorem]{Problem}
\newtheorem{question}[theorem]{Question}
\theoremstyle{definition}

\newcommand{\induce}[2]{\mbox{$ #1 \langle #2 \rangle$}}
\newcommand{\C}[1]{{\protect\cal #1}}

\begin{document}

\title{Problem collection from the IML programme: Graphs, Hypergraphs, and Computing}

\author{Klas Markstr\"om }
\maketitle

\begin{abstract}
	This is a collection of open problems and conjectures from the seminars and problem sessions of the 2014  
	IML programme: Graphs, Hypergraphs, and Computing.
\end{abstract}

\tableofcontents

\section{Introduction}
This collection of problems and conjectures is based on a subset of the open problems from the seminar series and the problem sessions of the IML programme Graphs, Hypergraphs, and Computing.  Each problem contributor has provided a write up of their proposed problem and the collection has been edited by Klas Markstr\"om.

\newpage
\section{Seminar  January 16, 2014, J\o{}rgen Bang-Jensen}

\subsection{Arc-disjoint spanning strong subdigraphs and disjoint Hamilton cycles}

\begin{conjecture}[Kelly 196?]
The arc set of every regular tournament can be decomposed into Hamilton cycles.
\end{conjecture}

\begin{theorem}[K\"uhn and Osthus, 2012]
The Kelly Conjecture is true for tournaments on $n$ vertices where $n\geq M$ for some very large $M$.
\end{theorem} 

As every $k$-regular tournament is 
$k$-arc-strong, the following Conjecture implies the Kelly conjecture.

\begin{conjecture}[Bang-Jensen and Yeo, 2004]\label{BJYTconj}
The arc set of every $k$-arc-strong tournament 
$T=(V,A)$ can be decomposed into $k$ disjoint sets
$A_1,\ldots{},A_k$ such that each of the spanning subdigraphs $D_i=(V,A_i)$, $i=1,2,\ldots{},k$
 is strongly connected.  
\end{conjecture}

\begin{theorem}[Bang-Jensen and Yeo 2004]\label{BJYthm}
Conjecture \ref{BJYTconj} is true in the following cases:
\begin{itemize}
\item when $k=2$
\item when every vertex of $T$ has in- and out-degree  at least $37k$
\item when there exists a non-trivial (both sides of the cut have at least 2 vertices) arc-cut of size $k$.
\end{itemize}
\end{theorem}

\begin{conjecture}[Bang-Jensen and Yeo, 2004]
There exists a natural number $K$ such that every $K$-arc-strong digraph $D=(V,A)$ can be decomposed 
into two arc-disjoint strong spanning subdigraphs $D_1=(V,A_1)$ and $D_2=(V,A_2)$.
\end{conjecture} 

The conjecture is true for tournaments with $K=2$ by Theorem \ref{BJYthm} and for locally semicomplete digraphs\footnote{a digraph is locally semicomplete if the out-neighbourhood and the in-neighbourhood of every vertex induces a semicomplete digraph. A digraph is semicomplete if it has no pair of non-adjacent vertices.}
with $K=3$ by a recent (non-trivial) result of Bang-Jensen and Huang (JCTB 2012).

For 2-regular digraphs (all in- and out-degrees equal to 2) the existence of arc-disjoint spanning
strong subdigraphs is equivalent to the existence of arc-disjoint Hamilton cycles. Hence, by the following theorem, it is NP-complete to decide whether a digraph has a pair of arc-disjoint strong spanning subdigraphs.
\begin{theorem}[Yeo, 2007]
It is NP-complete to decide whether a 2-regular digraph contains two arc-disjoint Hamilton cycles.
\end{theorem}

\begin{theorem}[K\"uhn, Osthus, Lapinskas and Patel, 2013]
There exists a natural number  $C$ such that every $Ck^2\log^2{}k$-strong tournament contains $k$ arc-disjoint Hamilton cycles. This is best possible up to the log-factor.
\end{theorem}

\begin{conjecture}[Thomassen, 1982]
Every 3-strong\footnote{A digraph $D$ is $k$-strong if it has at least $k+1$ vertices and $D-X$ is strongly connected for every subset $X$ of $V(D)$ of size at most $k-1$.} tournament has 2-arc-disjoint 
Hamilton cycles. 
\end{conjecture}

\subsection{Decompositions into vertex disjoint pieces/digraphs}

\begin{theorem}[K\"uhn, Osthus and Townsend, 2014]
There exists a naural number $C$ such that the vertex set of every $Ck^7t^4$-strong tournament
$T=(V,A)$ can be decomposed into disjoint subsets $V_1,V_2,\ldots{},V_t$  such that the tournaments
$T_i=\induce{T}{V_i}$ are $k$-strong for $i=1,2,\dots{},t$.
\end{theorem} 

\begin{question}
Can we also specify $t$ vertices $x_1,x_2,\ldots{},x_t$ and find $V_1,\ldots{},V_t$ as above such that
$x_i\in V_i$ holds for $i=1,2,\ldots{},t$?
\end{question}

\begin{conjecture}[Bermond and Thomassen, 1980]\label{BTconj}
Every digraph with minimum out-degree $2k-1$ contains $k$ disjoint directed cycles.
\end{conjecture}

For $k=2$ this was verified in 1983 by Thomassen who also proved the existence of a function $f(k)$ such that every digraph out minimum out-degree at least $f(k)$ has $k$ disjoint cycles.

The bound $2k-1$ is best possible as seen by considering the complete digraph on $2k-2$ vertices.

\begin{theorem}[Bang-Jensen, Bessy and Thomass\'e, 2013]
Conjecture \ref{BTconj} holds for tournaments
\end{theorem}

\subsection{Further open problems that were mentioned}

\begin{question}
Is there a polynomial algorithm for deciding whether the underlying graph $UG(D)$ of 
a digraph $D$ contains a 2-factor $C_1,C_2,\ldots{},C_k$ such that $C_1$ is a directed cycle in $D$,
while $C_i$, $i>1$ does not have to respect the orientations of arcs in $D$?
\end{question}

\begin{question}
What is the complexity of the following problem: given a 2-edge-coloured bipartite graph $B=(U,V,E)$; decide whether $B$ has two edge-disjoint perfect matchings $M_1,M_2$ so that every edge of $M_1$ has colour 1, while $M_2$ may use edges of both colours?
\end{question}

\newpage
\section{Seminar  January 16, 2014, Oleg Pikhurko}

Here are two open problems (as simplified as possible without losing their essence) that will be quite useful for measurable edge-colourings of graphings. --Oleg Pikhurko

\subsection{An Open Question about Finite Graphs}

\begin{problem} Estimate the minimum $f=f(d)$ such that the following holds. Let $G$ be a (finite) graph of maximum degree at most $d$ with at most $d$ pendant edges pre-coloured (with no two incident pre-coloured edges having the same colour). Then this pre-colouring can be extended to a proper  $(d+f)$-edge-colouring of $G$.\end{problem}

We can show $f=O(\sqrt d)$ suffices but it would be nice to prove that $f=O(1)$ is enough.

\subsection{Towards a Measurable Local Lemma}

For our purposes, it is enough to define a \emph{graphing} $\C G$ as a graph whose
vertex set is the unit interval $I=[0,1]$ (with the Borel $\sigma$-algebra $\C B$ and
the Lebesgue measure $\mu$) and whose edge set $E$ can be represented as
 $$
  E=\{\{x,y\}\mid x,y\in I,\ x\not=y,\ \exists i\in[k]\ \phi_i(x)=y\},
   $$
   for some (finite) family of measure-preserving invertible maps $\phi_1,\dots,
\phi_k:I\to I$. The general definition (and an excellent introduction) to graphings can
be found in Lov\'asz book \cite[Chapter 18]{lovasz:lngl}.

\begin{problem}\label{pr:1} Let $d\to\infty$. Prove that there is $g(d)=o(d/\log d)$ such that any graphing $\C G$ of maximum degree at most $d$ admits a measurable
partition $I=A\cup B$ such that for \emph{every} vertex $x$ is \emph{$(A,B)$-balanced}, meaning that its degrees into $A$ and into $B$ differ
by at most $g(d)$.\end{problem}

Some remarks: 
 \begin{enumerate}
 	\item The (finite) Local Lemma shows that the required partition $A\cup B$ exists for every finite graph
with $g(d)=O(\sqrt{d\log d})$. By the Compactness Principle, this extends to all countable graphs.

	 \item In Problem~\ref{pr:1} it is enough to find a partition such that the measure of the set $X$ of $(A,B)$-unbalanced vertices $x$ is zero. Indeed,  one can show that the measure of the union $Y$ of connectivity components of a graphing that intersect the null set $X$ is zero too. By the previous remark, we
can find a good partition of each component in $Y$; assuming the Axiom of Choice we can modify $A,B$
on the null set $Y$ to make \emph{every} vertex $(A,B)$-balanced.

	 \item Gabor Kun~\cite{kun:13:talk} proved some analytic version of the Local Lemma that in particular
implies for Problem~\ref{pr:1} that, for every $\varepsilon>0$, there is a measurable partition $I=A\cup B$ such that
the measure of $(A,B)$-unbalaned vertices is at most $\varepsilon$. But we do need the bad set to have measure zero in our application.
 \end{enumerate}

 \newpage
\section{Problem session Febuary 6, 2014}
 
 \subsection{Peter Allen}
 
Let $R(G,G)$ denote the 2-colour Ramsey-number for $G$.    It is known that there exists a constant $C$ such  that if $G$ is a planar graph on $n$ vertices then $R(G,G)\leq C n$.  It is also known that $C$ must be at least 4.

\begin{question} 
	Is $C \leq 12$?
\end{question}

\subsection{Hal Kierstead} 
An \emph{equitable} \emph{coloring} of a graph is a partiton of its
vertices into independent sets differing in size by at most one. In
1970 Hajnal\emph{ }and Szemerédi \cite{HSz} proved that for every
graph $G$ and integer $k$, if $\Delta(G)<k$ then $G$ has an equitable
$k$-coloring. Their proof did not yield a polynomial algorithm. About
$35$ years later, Mydlarz and Szemerédi, and independently Kostochka
and I, found such algorithms, and then joined forces to produce an
$O(kn^{2})$ algorithm \cite{KKMSz}. 

The \emph{maximum Ore degree} of a graph $G$ is $\theta(G):=\max\{d(x)+d(y):xy\in E(G)\}$.
Kostochka and I \cite{KK} proved that for every graph $G$ and integer
$k$, if $\theta(G)<2k$ then $G$ has an equitable $k$ coloring,
but the proof does not yield a polynomial algorithm.$ $ 
\begin{problem}
Is there a polynomial algorithm for constructing an equitable $k$-coloring
of any graph $G$ with $\theta(G)<2k?$
\end{problem}

\subsection{Jan van den Heuvel}

\subsubsection{Cyclic Orderings}

The following is a very special case of a much more general conjecture
which appears in \textsc{Kajitani et al.} (1988).

\begin{conjecture}\label{c1.1}\mbox{}\\*
  Let $T_1,T_2,T_3$ be edge-disjoint spanning trees in a graph~$G$ on~$n$
  vertices (so each tree has $n-1$ edges). Then there exists a cyclic
  ordering of the edges in $E(T_1)\cup E(T_2)\cup E(T_3)$ such that every
  $n-1$ cyclically consecutive edges in that ordering form a spanning tree.
\end{conjecture}

In fact, the same question can be asked for any number of spanning trees.
For two trees the result is proved in \textsc{Kajitani et al.} (1988), who
in fact prove it in the stronger form according to Conjecture~\ref{c1.2}
below.

Conjecture~\ref{c1.1} is really a problem about matroids. The following
appears in several places, including \textsc{Gabov} (1976),
\textsc{Cordovil~\& Moreira} (1993) and \textsc{Wiedemann} (2006).

\begin{conjecture}\label{c1.2}\mbox{}\\*
  Let $B=\{b_1,\dots,b_r\}$ and $B'=\{b'_1,\dots,b'_r\}$ be two disjoint
  bases of a matroid. Then there is a permutation
  $(b_{\pi(1)},\ldots,b_{\pi(r)})$ of the elements of~$B$ and a permutation
  $(b'_{\pi'(1)},\ldots,b'_{\pi'(r)})$ of the elements of~$B'$ such that
  the combined sequence
  $(b_{\pi(1)},\ldots,b_{\pi(r)},b'_{\pi'(1)},\ldots,b'_{\pi'(r)})$ is a
  cyclic ordering in which every~$r$ cyclically consecutive elements form a
  base.
\end{conjecture}

A weaker form of Conjecture~\ref{c1.2} is to just ask for a cyclic for a
cyclic ordering of $B_1\cup B_2$ (so we don't require that each base
appears as a consecutive part of the ordering). Even that conjecture is
open for matroids in general.

The most general conjecture in this area can be found in \textsc{Kajitani
  et al.} (1988); partial and related results appear in \textsc{van den
  Heuvel~\& Thomass\'e} (2012).

\bigskip
\noindent
\textsc{R. Cordovil and M.L. Moreira} \textsl{Bases-cobases graphs and
  polytopes of matroids}. Combinatorica~\textbf{13} (1993), 157--165.

\smallskip
\noindent
\textsc{H. Gabow}, \textsl{Decomposing symmetric exchanges in matroid
  bases}. Math.\ Programming~\textbf{10} (1976), 271--276.

\smallskip
\noindent
\textsc{J. van den Heuvel and S. Thomass\'e}, \textsl{Cyclic orderings and
  cyclic arboricity of matroids}. J.~Combin.\ Theory Ser.~B~\textbf{102}
(2012), 638-646.

\smallskip
\noindent
\textsc{Y. Kajitani, S. Ueno, and H. Miyano}, \textsl{Ordering of the
  elements of a matroid such that its consecutive~$w$ elements are
  independent}. Discrete Math.~\textbf{72} (1988) 187--194.

\smallskip
\noindent\textsc{D. Wiedemann}, \textsl{Cyclic base orders of
    matroids}. Manuscript, 2006. Retrieved 23~April 2007 from
  \url{http://www.plumbyte.com/cyclic_base_orders_1984.pdf}.\\
  Earlier version\,: \textsl{Cyclic ordering of matroids}. Unpublished
  manuscript, University of Waterloo, 1984

\subsubsection{Strong Colourings of Hypergraphs}

All hypergraphs in this section are allowed to have multiple edges and
edges of any size. The \textit{rank~$r(H)$} of a hypergraph~$H$ is the size
of the largest edge.

The \textit{strong chromatic number~$\chi_s(H)$} of a hypergraph~$H$ is the
smallest number of colours needed to colour the vertices so that for every
edge the vertices in that edge all receive a different colour. (So this is
the same as the chromatic number of the graph obtained by replacing every
edge by a clique.)

Let's call a \textit{derived graph of a hypergraph~$H$} a graph~$G$ on the
same vertex set as~$H$ where for each edge~$e$ of~$H$ of size at least two
we choose a pair $u,v\in e$ and add the edge~$uv$ to~$G$. And let's call
the following parameter the \textit{graph chromatic number of~$H$}:
\[\chi_d(H)=\max\{\,\chi(G)\mid\text{$G$ is a derived graph of
  $H$}\,\}.\]

It is obvious that $\chi_s(H)\ge\max\{r(H),\chi_d(H)\}$, and it is not so
hard to prove that $\chi_s(H\le\chi_d(H)^{\binom{r(H)}{2}}$. A little bit
more thinking will give
\[\chi_s(H)\le\chi_d(H)^{r(H)-1}.\]

The question is the find better upper bounds of~$\chi_s(H)$ in terms
of~$\chi_d(H)$ and~$r(H)$. It might even be true that there is an upper
bound that is linear in~$r(H)$.

\begin{question}\label{q2.1}\mbox{}\\*
  Does there exist a function $f:\mathbb{N}\rightarrow\mathbb{R}_+$ so
  that for every hypergraph~$H$ we have
  \[\chi_s(H)\le f(\chi_d(H))\cdot r(H)\,\text{?}\]  
\end{question}

A neat argument, due to my former PhD student Alexey Pokrovskiy, gives a
proof that if $\chi_d(H)=2$, then $\chi_s(H)=r(H)$.

Other special classes of hypergraphs for which Question~\ref{q2.1} has a
positive answer can be found in \textsc{Dvo\v{r}\'ak~\& Esperet} (2013).
That paper was also the inspiration for starting to think about this type
of questions.

\bigskip
\noindent
\textsc{Z. Dvo\v{r}\'ak and L. Esperet}, \textsl{Distance-two coloring of
  sparse graphs}. arXiv:1303.3191 [math.CO],
\url{http://arxiv.org/abs/1303.3191} (2013), 13~pages.

\subsection{Victor Falgas-Ravry} 
\subsubsection{Largest antichain in the independence complex of a graph}

Write $Q_n$ for the collection of all subsets of $[n]=\{1,2, \ldots n\}$. We denote by $Q_n^{(r)}$ the \emph{$r^{\textrm{th}}$ layer} of $Q_n$, that is, the collection of all subsets of $[n]$ of size $r$.

A family $\mathcal{A}\subseteq Q_n$ is an \emph{antichain} if for every pair of distinct elements $A,B \in \mathcal{A}$, $A$ is not a subset of $B$ and $B$ is not a subset of $A$. How large an antichain can we find in $Q_n$? Clearly each layer of $Q_n$ forms an antichain, and a celebrated theorem of Sperner from 1928 asserts that we cannot do better than picking a largest layer:
\begin{theorem}[Sperner~\cite{Sperner28}]
Let $\mathcal{A}\subseteq Q_n$ be an antichain. Then 
\[\vert  \mathcal{A} \vert \leq \max_r \vert Q_n^{(r)}\vert.\]
\end{theorem}

I am interested in a generalisation of Sperner's theorem where $\mathcal{A}$ is restricted to a subset of $Q_n$: suppose we are given a graph $G$ on $[n]$. The \emph{independence complex of $G$} is the collection of all independent sets from $V(G)=[n]$,
\[Q(G)=\{A \subseteq [n]: \ A \textrm{ independent in }G\}.\]
We write $Q^{(r)}(G)= \{A \in Q(G): \ \vert A\vert =r\}$ for the \emph{$r^{\textrm{th}}$ layer} of $Q(G)$, and define the \emph{width} of $G$, $s(G)$, to be the size of a largest antichain in $Q(G)$. Clearly we have
\begin{align}
s(G)&\geq \max_r \vert Q^{(r)}(G) \vert.\label{trivialbound}
\end{align}
In general, $s(G)$ can be much larger than this: it is not hard to construct examples of graph sequences $(G_n)_{n\in \mathbb{N}}$ for which $\max_r \vert Q^{(r)}(G_n\vert) =o\left(s(G_n)\right)$. However one would expect that if $G$ is reasonably homogeneous then (~\ref{trivialbound}) should be close to tight.
\begin{question}
When do we have (almost) equality in (~\ref{trivialbound}) ? What conditions on $G$ are sufficient to guarantee (almost) equality ?
\end{question}

I am particularly interested in the cases where $G=C_n$, the cycle of length $n$, or where $G=P_n$, the path of length $n-1$. In this setting, an analogue of the Erd\H{os}--Ko--Rado theorem in $Q(G)$ was proved by Talbot~\cite{Talbot03}, using an ingenious compression argument. It is known~\cite{FalgasRavry14} that the size of a largest antichain in a class of graphs including both $C_n$ and  $P_n$ is of the same order as the size of a largest layer. However we really should have equality here:

\begin{conjecture}
\[s(C_n)= \max_r \vert Q(C_n)^{(r)} \vert \textrm{ and }s(P_n)= \max_r \vert Q(P_n)^{(r)}\vert.\]
\end{conjecture}

It would also be interesting to know what happens in the case of random graphs:
\begin{question}
Let $p=cn^{-1}$, for some constant $c>0$. Is it true with high probability that
\[s(G_{n,p})=(1+o(1)) \max_r \vert Q(G_{n,p})^{(r)} \vert? \]
\end{question}

\newpage
\section{Problem session Febuary 19, 2014}

\subsection{Mikl\'os Simonovits}

\begin{problem}[Paul Erd\H{o}s, via M. Simonovits]
  Is it true that if $G_n$ is a 4-chromatic $n$-vertex graph and 
  deleting any edge of it we get a 3-chromatic graph, then the
  minimum degree of $G_n$ is $o(n)$ (as $n\to\infty$)?
\end{problem}

{\bf Motivation, partial results: } A graph $G$ is called {\bf
  $k$-color-critical} if it is $k$-chromatic but deleting any edge of
$G$ we get a $k-1$-chromatic graph. (Actually, we could speak of
edge-critical and vertex-critical graphs, but  we stick to the
edge-critical case.)

Bjarne Toft and myself, using a construction of Toft and a
transformation of mine constructed (infinitely many) 4-colour-critical graphs where the
minimum degree is $>c\root 3 \of n$.  (Simonovits, M.: On
colour-critical graphs. Studia Sci. Math. Hungar. 7 (1972), 67--81, and
Toft, B. Two theorems on critical 4-chromatic graphs, Studia
Sci. Math. Hungar. 7 (1972), 83--89.) I do not know of anything with
higher minimum degree (though I may overlook some newer results?)

A trivial construction of G. Dirac, obtained by joining two odd
$n/2$-cycles completely shows that there exist $6$-critical graphs
with minimum degrees $n/2+2$. The difficulties occur for 4 and
5-critical graphs.  (The 3-critical graphs are just the odd cycles.)

As I wrote, a basic ingredient of our construction was an earlier
construction of Bjarne Toft, a 4-critical graph with $\approx
{n^2\over 16}$ edges, where the vertices are in four groups of $n/4$
vertices, and (a) the first and last groups form two odd cycles, (b) the
second and third groups form a complete bipartite graphs, (c) the first
group is joined to the second class by a 1-factor and the third group
to the last group also by a 1-factor.

\subsection{Fedor Fomin} 

\begin{question}
For a given $n$-vertex planar graph $G$, is it possible to find in polynomial time (or to show that this is NP-hard) an independent set of size $\lfloor n/4 \rfloor +1$?  
\end{question}

\subsection{Carsten Thomassen} 

Smith's theorem says that, for every edge e in a cubic graph, there is in an even number of Hamiltonian cycles containing e. As a consequence, every cubic Hamiltonian graph has at least 3 Hamiltonian cycles.
 
If e is an edge of a Hamiltonian bipartite cubic graph G, then G has en even number of Hamiltonian cycles through e. Using Smith's theorem once more, also G-e han an even number of Hamiltonian cycles. Hence the total number of Hamiltonian cycles is even, in contrast to the situation for non-bipartite graphs where it is 3 for infinitely many graphs.

\begin{problem} 
Does there exist a 3-connected cubic bipartite graph having an edge e such that there are precisely two Hamiltonian cycles containing e? (Otherwise, there will be at least 4 such cycles.)
 \end{problem}
 
 \begin{problem} 
 Does there exist a 3-connected cubic bipartite graph having precisely 4 Hamiltonian cycles? (Otherwise, there will be at least 6 such cycles.)
 \end{problem}
 
It is an old problem whether a second Hamiltonian cycle in a Hamiltonian cubic graph can be found in polynomial time. Andrew Thomason''s lollipop method is a simple algorithm producing a second Hamiltonian cycle, but it may take exponential many steps. The known examples have many edge-cuts with three edges.
 
  \begin{problem} \label{thp}
 Does there exist a family of cubic, cyclically 4-edge-connected Hamiltonian graphs for which the lollipop method takes superpolynomially many steps?
  \end{problem}
  
If Problem \ref{thp} has a negative answer, one can show that there exists a polynomially bounded algorithm for finding a second Hamiiltonian cycle in a cubic Hamiltonian graph.

\subsection{J\"orgen Backelin}
The problem, stated briefly:

Determine the exact chromatic numbers for shift graphs with short vertices for cyclically ordered points.

Consider a finite set $C$, ordered cyclically. Define a shift graph by letting the vertices be all
$r$-subsets of $C$, for some small $r$ (I suggest doingt his for $r=2$ or 3, in the first place), and by
letting two vertices form an edge if their are disjoint, and intertwined in a prescribed manner.
\begin{problem}
	Determine the exact chromatic numbers for these graphs.
\end{problem}

Detailed explanation:

A cyclic order on a finite set
$C$ intuitively is what you think it should be, if you place the elements in a circle and follow it e. g.
clockwise: It is not meaningful to say that the element a precedes b; but it is meaningful to say
that starting from $a$ we pass $b$ before encountering $c$. Thus, a formal definition would have to deal with ternary rather than binary predicates (properties).

Cyclic orders was treated by P. J. Cameron 1978 (Math. Z. 148, pp. 127-139). One way to define them
formally is as a property P, which holds for some triples of different elements in C, such that for
any different $a, b, c, d$ in $C$ we have:
\begin{itemize}
	\item	 Precisely one of $P(a,b,c)$ and $P(a,c,b)$ holds.
	\item If $P(a,b,c)$, then $P(b,c,a)$.
	\item If $P(a,b,c)$ and $P(a,c,d)$, then $P(a,b,d)$.
\end{itemize}
(Think of $P(a,b,c)$ as the statement "Starting from $a$, we pass $b$ before arriving at $c$".)

"Intertwining vertices" should be done analogously as for ordinary shift graphs. E. g., for $r = 2$, there are two possibilities for two disjoint
vertices $\{a,b\}$ and $\{c,d\}$: Either the pattern XXOO, which holds if either both $P(a,b,c)$ and $P(a,b,d)$ hold, or neither does, or the pattern XOXO. The
first case has a trivial chromatic number for ordinary shift graphs, but not self-evidently so in the cyclic order situation.

For $r=3$, there are essentially only three intertwinement patterns: XXXOOO, XXOXOO, and XOXOXO.

Rationale:

I do not know if there are any direct applications of this. In general, working with cyclic orders removes a kind of lack of balance between different parts 
of a graph, which means that we often may construct more efficient examples of graphs with certain properties in this manner. Thus, I find this a more natural setting.

Note on generalisations:
There are cyclically ordered sets of any cardinality. I do not know if there is a theory for "cyclically well-ordered sets" of large cardinalities, though.

\subsection{Andrzej Ruci\'nski} 

For graphs $F$ and $G$ and a positive integer $r$, we write $F\to(G)_r^v$ if every $r$-coloring of the vertices of $F$ results in a monochromatic copy of $G$ in $F$ (not necessarily induced).  E.g., $K_{r(s-1)+1}\to(K_s)^v_r$ by the Pigeon-hole Principle, where $K_n$ is the complete graph on $n$ vertices. Define
$$mad(F)=\max_{H\subseteq F}\frac{2|E(H)|}{|V(H)|}$$
and 
$$m_{cr}(G,r)=\inf\{mad(F)\;:\;F\to (G)_r^v\}.$$
\begin{problem} 
	Determine or estimate $m_{cr}(G,r)$ for every graph $G$ and $r\ge2$.
\end{problem}

\begin{figure}[!t]
	\begin{center}
		\includegraphics[width=1\textwidth]{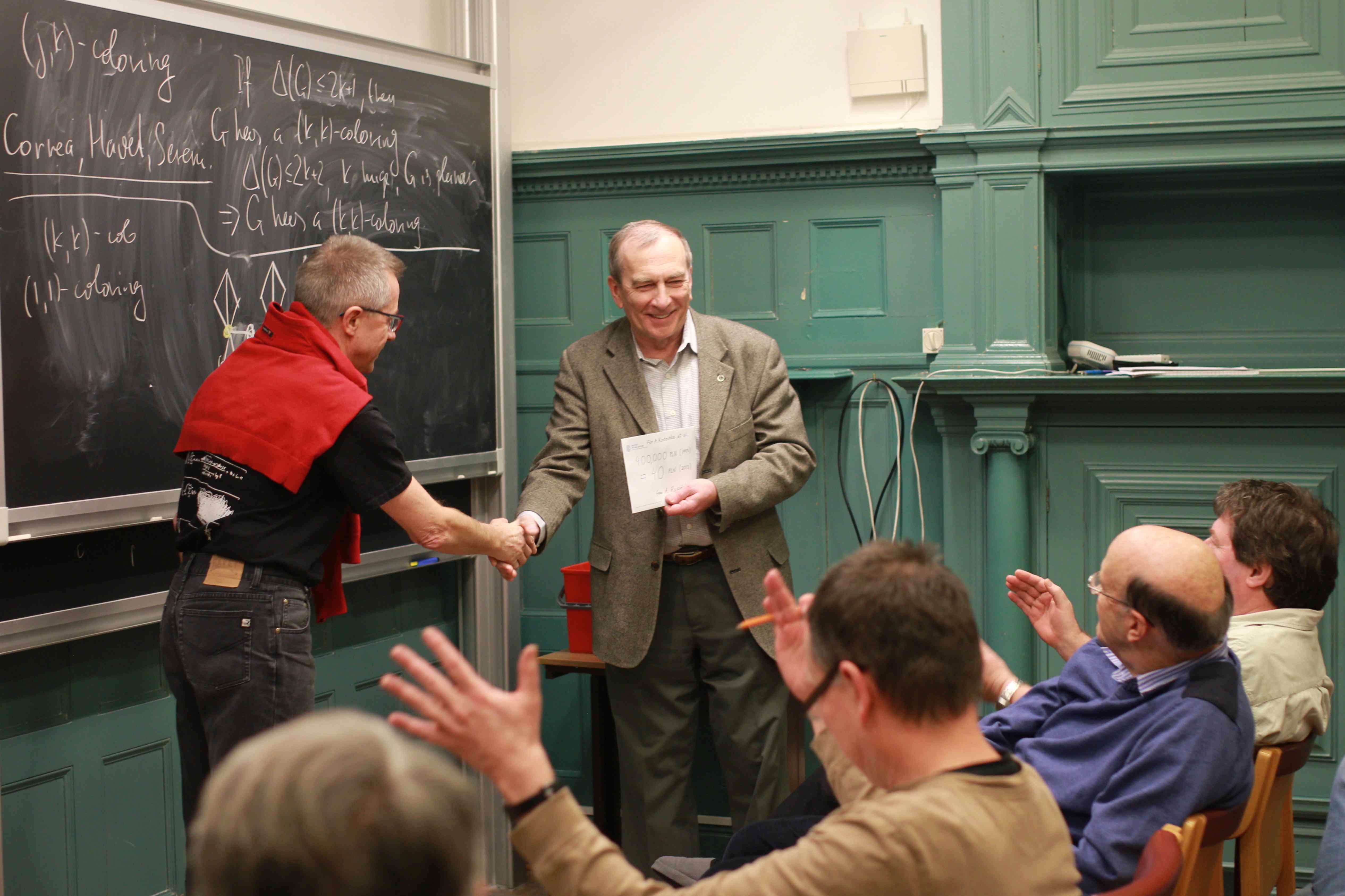}
		\caption{Sasha Kostochka receiving the award}
		\label{fig4}
	\end{center}
 \end{figure} 

It is known \cite{kr} that
$$r\max_{H\subseteq G}\delta(H)\le m_{cr}(G,r)\le 2r\max_{H\subseteq G}\delta(H),$$
where $\delta(G)$ is the minimum vertex degree in $G$.
The lower bound is attained by complete graphs, that is, $m_{cr}(K_s,r)=mad(K_{r(s-1)+1})=r(s-1)$.
On the other hand, the upper bound is asymptotically achieved by large stars as it was proved in \cite{kr} that, in particular for $r=2$ colors, 
$$4-\frac4{k+1}\le m_{cr}(S_k,2) \le 4-\frac{2(k+1)}{k^2+1},$$
where $S_k$ is the star with $k$ edges.
For $k=2$ this reads
$$\frac83\le m_{cr}(S_2,2) \le\frac{14}5.$$
In \cite{kr} I offered 400,000 z³ (Polish currency in 1993, equivalent after denomination of 1995 to 40 PLN) for the determination of the exact value of
$ m_{cr}(S_2,2)$.
Recently, it was pointed out by A.Pokrovskiy that a result of  Borodin, Kostochka, and Yancey \cite{sasha} 1-improper 2-colorings of sparse graphs yields that
$$m_{cr}(S_2,2)=\frac{14}5.$$
During the Open Problem session at the Mittag-Leffler Institute I handed to Sasha Kostochka an envelope with four 10 PLN bills so that he could share the award among his co-authors and A. Pokrovskiy as well.

\begin{problem}  
	Determine $m_cr(S_2,3)$ (no monetary award for that!) The current bounds, unimproved  since 1994 are
	$$\frac{18}5\le m_cr(S_2,3)\le \frac{22}5.$$
\end{problem}

\newpage
\section{Seminar  February 20, 2014,  Allan Lo}

For a graph $G$, we know that $\chi'(G) = \Delta (G)$ or $\chi'(G) = \Delta (G)+1$.
A subgraph $H$ of $G$ is \emph{overfull} if $e(H) > \Delta (G) \lfloor |H|/2 \rfloor$ (note this requires $|H|$ to be odd). 
Note that an overfull subgraph is a trivial obstruction for $\chi'(G) = \Delta(G)$.
In 1986, Chetwynd and Hilton~\cite{overfull} gave the following conjecture.
\medskip

\begin{conjecture}[Overfull subgraph conjecture]
A graph $G$ on $n$ vertices with $\Delta (G) \geq n/3$ satisfies $\chi'(G)=\Delta(G)$ if and only if $G$ contains no overfull subgraph.
\end{conjecture}

(Some remark about regular graphs)
The $1$-factorization conjecture is a special case of the overfull subgraph conjecture.
If $G$ is $d$-regular and contains no overfull subgraph, then ($|G|$ is even) and every odd cut has size at least $d$ edges.
So $G$ has a $1$-factor.
Meredith~\cite{meredith} showed that for all $d\ge 3$, there exists a $d$-regular graph $G$ on $20d-10$ vertices with $\chi'(G) = d+1$, which contains no overfull subgraph.

 \newpage
\section{Problem session March 5, 2014}

\subsection{Jacques Verstra\"ete} 

\subsubsection*{A problem on Majority Percolation}

Let $G$ be a finite graph, and let $p \in [0,1]$. Suppose that vertices of $G$ are randomly and independently {\em infected} with probability $p$ -- this is the {\em infection probability}. Then consider the following deterministic rule: at any stage, an uninfected vertex becomes infected if {\em strictly more} than half of its neighbors are
infected. Let $A(G)$ be the event that in finite time every vertex of $G$ becomes infected with associated probability measure $P_p$.

\begin{problem} 
	Does there exist a sequence of graphs $(G_n)_{n \in \mathbb N}$ such that for every $p > 0$:
	\[ \lim_{n \rightarrow \infty} P_p(A(G_n)) = 1 \quad \mbox{?}\]

\end{problem}

I believe the answer is no. This process is a version of a process called {\em bootstrap percolation}. 
The most studied case is the $n \times n$ grid $\Gamma_n$, with the rule that {\em at least} two infected neighbors of an uninfected vertex cause the vertex to become infected. In this case, a remarkable paper of Holroyd shows that for all $\varepsilon > 0$, 
\[ P_p(A(\Gamma_n)) \rightarrow \left\{\begin{array}{ll}
1 & \mbox{ if } p > (1 + \varepsilon)\frac{\pi^2}{18\log n}  \\
0 & \mbox{ if } p < (1 - \varepsilon)\frac{\pi^2}{18\log n}
\end{array}\right.
\]
Finer control of the relationship between $\varepsilon$ and $n$ was obtained by Graver, Holroyd and Morris.

\subsection{Klas Markstr\"om} 

Given a matrix $m \in \mathrm{GL}(n,2)$, i.e an invertible matrix with entries 0/1, we let $\mathcal{D}(m)$ denote the smallest number of row-operations we can use in order to reduce $m$ to the identity matrix. Equivalently  $\mathcal{D}(m)$ is the distance from $m$ to $I$ in the Cayely graph $\mathrm{Cay}(\mathrm{GL}(n,2),S)$, where $S$ is the set of elementary matrices.  In \cite{AHM}  an algorithm was given which can row reduce a matrix $m$ to te identity using 
$$\frac{n^2}{\log_{2} n} + o\left(\frac{n^2}{\log_{2} n}\right ) $$
row operations, and it was proven that the expected value of $\mathcal{D}(m)$ for a random matrix from $\mathrm{GL}(n,2)$ is not less than half of that. Equivalently, this shows that the diameter of $\mathrm{Cay}(\mathrm{GL}(n,2),S)$ is at most the first bound, and the average distance is at least the second.

\begin{problem}{\ }\\
	\begin{enumerate}
		\item Give an explicit (non-random) example of an matrix  $m \in \mathrm{GL}(n,2)$ such that $\mathcal{D}(m)\geq 100 n$
		
		\item Give an explicit  example of an matrix  $m \in \mathrm{GL}(n,2)$ such that $\mathcal{D}(m)\geq  n\log n$
		
	\end{enumerate}
\end{problem}

\subsection{Andrzej Ruci\'nski} 

Given integers $1\leq \ell< k$, we define an {\em $\ell$-overlapping cycle} as a $k$-uniform
hypergraph (or $k$-graph, for short) in which,  for some cyclic ordering of its vertices, every
edge consists of $k$ consecutive vertices, and every two consecutive edges (in the natural ordering
of the edges induced by the ordering of the vertices) share exactly $\ell$ vertices. If $H$
contains an $\ell$-overlapping Hamiltonian cycle then $H$ itself is called
\emph{$\ell$-Hamiltonian}.

 A $k$-graph $H$ is
\emph{$\ell$-Hamiltonian saturated}, $1\le \ell\le k-1$, if $H$ is not $\ell$-Hamiltonian  but for
every $e\in H^c$ the $k$-graph $H+e$ is such. For $n$ divisible by $k-\ell$, let $sat(n,k,\ell)$ be
the \emph{smallest} number of edges in an $\ell$-Hamiltonian saturated $k$-graph on $n$ vertices.
In the case of graphs, Clark and Entringer \cite{CE} proved in 1983 that $sat(n,2,1)=\lceil
\tfrac{3n}2\rceil$ for $n$ large enough.

A. \.Zak showed  that for $k\ge2$, $sat(n,k,k-1)=\Theta(n^{k-1})$ \cite{zak}. Together, we proved
that for all $k\ge3$ and $\ell=1$, as well as for all  $\tfrac45k\le\ell\le k-1$
\begin{equation}\label{1}
sat(n,k,\ell)=\Theta(n^{\ell}),
\end{equation}
 and  conjectured that (\ref{1}) holds for all $k$ and $1\le\ell\le k-1$ \cite{RZ}. The smallest
open case is $k=4\,, \ell=2$. Recently, we have got some partial results: $sat(n,k,\ell)=O(n^{
(k+\ell)/2 })$ and $sat(n,4,2)=O( n^{14/5} )$.

 \newpage
\section{Problem session March 19, 2014}

\subsection{Brendan McKay} 
\noindent\textbf{How many $n/2$-cycles can a cubic graph have?}

Given a simple cubic graph with $n$ vertices, what is a good upper bound on the number of cycles of length $n/2$ it can have?

A random cubic graph has $\Theta((4/3)n/n)$ cycles of length $n/2$. So do random cubic bipartite graphs. Also the whole cycle space has
size $2^{n/2+1}$, so twice that is a (silly) upper bound.

The actual maximums for $n=4,6,\ldots,24$ are: 0,2,6,12,20,20,48,48,132,118,312 (not in OEIS). All these are achieved uniquely except that for 20 vertices
there are two graphs with 132 10-cycles.

\noindent\textbf{Maximum automorphism group for a 3-connected cubic graph}

Let $a(n)$ be the greatest order of the automorphism group of a
3-connected cubic graph with $n$ vertices. I conjecture: for $n\ge 16$, $a(n)<n\, 2^{n/4}$.

There is a paper of Opstall and Veliche that finds the maximum over all
cubic graphs, but the maximum occurs for graphs very far from being 3-connected.

When $n$ is a multiple of 4 there is a vertex-transitive cubic graph
achieving half the conjectured bound, so if true the bound is pretty sharp.

A paper of Poto\v cnik, Spiga and Verret (arxiv.org/abs/1010.2546), together with some computation, establishes the conjecture for vertex-transitive graphs, so the remaining problem is whether one can do better for non-transitive graphs. For 20, and all odd multiples of 2 vertices from 18 to at least 998 (but not for 4--16 or 24 vertices) the graph achieving the maximum is not vertex-transitive.

\noindent\textbf{Probability that a random integer matrix is positive}

Let $M(n,k)$ be the set of $n\times n$ matrices of nonnegative integers such that every row and every column sums to $k$. Let $P(n,k)$ be the fraction of such matrices which have no zero entries, equivalently the probability that a random matrix from the uniform distribution on $M(n,k)$ has no zero entries.

Obviously $P(n,k)=0$ for $k<n$.  It seems obvious that $P(n,k)$ should be increasing as a function of $k$ when $n$ is fixed and $k\ge n$, but can you prove it?

We know that $P(n,k)=|M(n,k)|/|M(n,k)|$.  Also, note that $M(n,k)$ is the set of integer points in the $k$-dilated Birkhoff polytope, and $P(n,k)$
is the fraction of such points that don't lie on the boundary.  Ehrhart theory tells us that $|M(n,k)|=H_n(k)$ where $H_n$ is a polynomial, and that
$P(n,k) = (-1)^{n+1} H_n(-k)/H_n(k)$.  Does it help?

\subsection{Andrew Thomason (by proxy)} 

Suppose you are given a set $E$ and a collection of finite sequences of
elements of $E$.  We now wish to determine if there is a graph such that $E$ is the
edges of the graph and the sequences are (nice, simple) paths in the graph.

This can be done, the graph would have at most $2|E|$ vertices so you
can search all possibilities. So the question would be whether you can do
it efficiently.

\subsection{Bruce Richter}
Haj\'os conjectured that if the chromatic number $\chi(G)$ of a graph $G$ is at least $r$, then $G$ contains a subdivision of $K_r$.  Hadwiger conjectured that $G$ has $K_r$ as a minor.

\bigskip{\bf Albertson's Conjecture}:  {\em If $\chi(G)\ge r$, then the crossing number cr$(G)\ge \textrm{cr}(K_r)$.}

\bigskip
This is known for $r\le 16$, with the best result being that of J. Bar\'at and G. T\'oth, Towards the Albertson conjecture, Elec.\ J.\ Combin.\ 2010.   The interesting thing here is that the crossing number of $K_r$ is only known when $r\le 12$.  

The conjecture is easy for large $r$-critical graphs.  The problem occurs when $|V(G)|$ is just a little larger than $r$.

\subsection{Alexander Kostochka}

\subsubsection{Problem 1}
For nonnegative integers $j,k$, a $(j,k)$-coloring of graph $G$ is a partition $V(G)=J\cup K$ such that
$\Delta(G[J])\leq j$ and $\Delta(G[K])\leq k$. An old result of Lov\' asz implies that every graph $G$ with
$\Delta(G)\leq j+k+1$ has a $(j,k)$-coloring. The proof is short and one may wonder about possible
Brooks-type refinements of the result. In seeking such possibilities, Corr\' ea,  Havet, and Sereni~[1],
conjectured that for sufficiently large $k$ (say, $k>10^6$), every planar graph $G$ with $\Delta(G)\leq 2k+2$ has a $(k,k)$-coloring.

\bigskip
[1] R. Corr\' ea,  F. Havet, and J.-S. Sereni, 
About a Brooks-type theorem for improper colouring. 
Australas. J. Combin. 43 (2009), 219--230.

\subsubsection{Problem 2}
A graph is a {\em circle graph}, if it is the
intersection graph of a family of chords of a circle. Circle graphs arise
in many combinatorial problems ranging from sorting problems to studying planar
graphs to continous fractions. In particular,
for a given permutation $P$ of $\{1,2,\ldots,n\}$,
the problem of finding the minimum number of stacks needed to
obtain the permutation $\{1,2,\ldots,n\}$  from $P$ reduces to finding
the chromatic number of a corresponding circle graph.
There are polynomial algorithms for finding the clique number and the independence number of a circle graph,
but finding the chromatic number of a circle graph is an NP-hard problem.

Let $f(k)$ denote the maximum chromatic number of a circle graph with clique number $k$.
Gy\' arf\' as~\cite{G2}
proved that $f(k)$ is well defined and $ f(k)\leq 2^k(2^k-2)k^2$. The only known exact
value is $f(2)=5$. The best bounds known to me are $$0.5k(\ln k-2)\leq f(k)\leq 50\cdot 2^k.$$
The lower bound is only barely superlinear, and the upper is very superlinear. It would be
interesting to improve any of them. More info on circle graphs and their colorings could
be found in~\cite{Gol,G1,AK}.

\subsection{David Conlon}   
Monochromatic cycle partitions in mean colourings

A well-known result of Erd\H{o}s, Gy\'arf\'as and Pyber~\cite{EGP91} says that there exists a constant $c(r)$, depending only on $r$, such that if the edges of the complete graph $K_n$ are coloured with $r$ colours, then the vertex set of $K_n$ may be partitioned into at most $c(r)$ disjoint monochromatic cycles, where we allow the empty set, single vertices and edges to be cycles. For $r = 2$, it is known \cite{BT10} that two disjoint monochromatic cycles of different colours suffice, while the best known general bound \cite{GRSS06} is $c(r) = O(r \log r)$.

With Maya Stein~\cite{CS14}, we recently considered a generalisation of this monochromatic cycle partition question to graphs with locally bounded colourings. We say that an edge colouring of a graph is an {\it $r$-local colouring} if the edges incident to any vertex are coloured with at most $r$ colours. Note that we do not restrict the total number of colours. Somewhat surprisingly, we prove that even for local colourings, a variant of the Erd\H{o}s-Gy\'arf\'as-Pyber result holds. 

\begin{theorem} \label{thm:rcolours}
The vertex set of any $r$-locally coloured complete graph may be partitioned into $O(r^2 \log r)$ disjoint monochromatic cycles.
\end{theorem}

For $r = 2$, we have the following more precise theorem.

\begin{theorem} \label{thm:2colours}
The vertex set of any $2$-locally coloured complete graph may be partitioned into two disjoint monochromatic cycles of different colours.
\end{theorem}

An edge colouring of a graph is said to be an $r$-mean colouring if the average number of colours incident to any vertex is at most $r$. We suspect that a theorem analogous to Theorem~\ref{thm:rcolours} may also hold for $r$-mean colourings but have been unable to resolve this question in general. 

\begin{question}
Does there exist a constant $m(r)$, depending only on $r$, such that the vertex set of any $r$-mean coloured graph may be partitioned into at most $m(r)$ cycles?
\end{question}

We can show that the vertex set of any $2$-mean coloured graph may be partitioned into at most two cycles of different colours but the proof uses tricks which are specific to the case $r = 2$.

 \newpage
\section{Problem session April 2nd, 2014}

\subsection{Mikl\'os Simonovits} 

The problem we discuss here is informally as follows:

Is it true that the Tur\'an number of an infinite family of forbidden
(bipartite) graphs $\mathcal{L}$ can be approximated arbitrarily well in the
exponent by finite subfamilies?

More precisely, here we consider ordinary simple graphs: no loops or
multiple edges are allowed.  ${\bf ext}(n,\mathcal{L})$ denotes the maximum number
of edges a graph $G_n$ on $n$ vertices can have without containing
subgraphs from $\mathcal{L}$.  The problem below is motivated by the fact that
if $\mathcal{C}$ is the family of all cycles then ${\bf ext}(n,\mathcal{C})=n-1$, however
for any finite $\mathcal{C}^*\subset\mathcal{C}$, for some $\alpha=\alpha(\mathcal{C}^*)>0$,
and $c_1>0$, ${\bf ext}(n,\mathcal{C}^*)>c_1n^{1+\alpha}$.  This means that some
kind of compactness is missing here. On the other hand, the continuity
of the exponent easily follows from Bondy-Simonovits theorem:
$${\bf ext}(n,C_{2k})=O(n^{1+{1/k}}).$$
This means a continuity in the exponent.

\begin{problem}[Continuity of the exponent]
Let $\mathcal{L}$ denote an infinite family of bipartite (excluded) graphs and 
$$\mathcal{L}_m:=\{L~:~L\in\mathcal{L},~ v(L)\le m\}.$$
Is it true that if for some $\alpha>0$ and $c>0$ we have ${\bf ext}t(n,\mathcal{L})=
O(n^{1+\alpha})$, then for any $\varepsilon>0$, we have
$${\bf ext}(n,\mathcal{L}_m)=O(n^{1+\alpha+\varepsilon}),\quad\mbox{as}\quad 
n\to\infty,$$
if $m$ is large enough?

\end{problem}

\subsection{Erik Aas and Brendan McKay}
Let $G = (V,E)$ be a connected simple graph, $k^E$ its edge space over the field $k$. We are interested in the subspace $C(G)$ spanned by the (characteristic vectors of the) cycles of $G$.

When $k$ is the field with two elements, it is a classical fact that the dimension of $C(G)$ is $|E(G)| - |V(G)| + 1$. This can be proved by providing an explicit basis of $C(G)$, as follows. Pick any spanning tree $T$ of $G$, and for each edge $e$ not in $T$ consider the unique cycle whose only edge not in $T$ is $e$. These cycles are linearly independent and thus span $C(G)$ in this case.

Now, when $k$ does not have characteristic $2$, the dimension is not a simple function of $|E(G)|$ and $|V(G)|$. However, in the case $G$ is $3$-edge-connected, it is not difficult to prove that in fact the dimension of $C(G)$ is $|E(G)|$.

{\bf Question: Is there a nice explicit basis for $C(G)$ consisting of cycles indexed by $E$, assuming $G$ is $3$-edge-connected?}

\subsection{Andr\'as Gy\'arf\'as}

 A {\em 3-tournament $T_n^3$ } is the set of all
triples on vertex set $[n]=\{1,2,\dots,n\}$ such that in each triple
some vertex is designated as the {\em root} of the triple. A set
$X\subset [n]$ is a {\em dominating set} in a 3-tournament $T_n^3$
if for every $z\in [n]\setminus X$ there exist $x\in X, y\in [n]$
($y\ne z, y\ne x$) such that $x$ is the root of the triple
$(x,y,z)$. Let $dom(T_n^3)$ denote the cardinality of a smallest
dominating set of $T_n^3$.

\begin{conjecture} 
	There exists a 3-tournament $T_n^3$ such that $dom(T_n^3)\ge 2014$.
\end{conjecture}

\begin{conjecture}
	 If any four vertices of a 3-tournament $T_n^3$ contain at least two triples with the same root then $dom(T_n)\le 2014$.
\end{conjecture}

I already posed this pair of conjectures at the 2012 Prague
Midsummer Combinatorial Workshop (of course with 2012 in the role of
2014).

Note that the $2$-dimensional versions of the above conjectures are
true: there exist tournaments $T$ with  $dom(T)\ge 2014$; if any
three vertices of a tournament $T$ contain two pairs with the same
root then $dom(T)=1$. Also, if {\em three triples} are required with
the same root in every four vertices of a $T_n^3$ then
$dom(T_n^3)=1$ follows easily (a remark with Tuza).

\subsection{Klas Markstr\"om}

If $G=(V,E)$ is an $n$-vertex graph then the strong chromatic number of $G$, denoted $s_\chi (G)$, is the minimum $k$ such that the
following hold: Any graph which is the union of $G$ and a set of $\lceil \frac{n}{k} \rceil$ vertex-disjoint $k$-cliques is $k$-colourable. Here we take the union
of edge sets, adding isolated vertices to G if necessary to make $n$ divisible by $k$. 

It is easy to see that $\Delta{G}+1 \leq s_\chi (G)$, and Penny Haxell \cite{P1} has proven that $s_\chi (G) \leq c \Delta(G)$ for all $c>\frac{11}{4}$ if $\Delta$ is large enough, and \cite{P2}  $s_\chi (G) \leq 3 \Delta(G) -1$ in general.   The folklore conjecture here is that  $s_\chi (G) \leq 2 \Delta(G)$.  This is known to be true if $\Delta \geq \frac{n}{6}$ \cite{JJM}.

Let us define the  \emph{biclique number} $\omega_b(G)$ to be the maximum $t$ such that there exists a $K_{a,b}\subset G$ with $t=a+b$.

A few years ago I made the following conjecture:
\begin{conjecture}
	$\omega_b(G)\leq s_\chi (G) \leq \omega_b(G) +1 $
\end{conjecture}
The lower bound is easily seen to be true so the conjecture really concerns the upper bound.

\subsection{Benny Sudakov}
\begin{question}
How many edges do we need to delete to make a $K_r$-free graph $G$ of 
order $n$ bipartite? 
\end{question}
For $r=3, 4$ this was asked long time ago by P. Erd\H{o}s. For triangle-free graphs he conjectured that deletion of 
$n^2/25$ edges is always enough and that extremal example is a blow-up of a $5$-cycle. Sudakov answered the question for $r=4$ and proved 
that the unique extremal construction in this case is a complete $3$-partite graph with equal parts. This result suggests that a complete 
$(r-1)$-partite graph of order $n$ with equal parts is worst example also for all remaining values of $r$. Therefore we believe that it is enough to 
delete at most $\frac{(r-2)^2}{4(r-1)^2}n^2$ edges for even $r\geq 5$ and at most $\frac{r-3}{4(r-1)}n^2$ edges for odd $r\geq 5$ to make bipartite 
any $K_r$-free graph $G$ of order $n$.

 \newpage
\section{Problem session April 16th, 2014}

\subsection{Dhruv Mubayi} 

Fix $k \ge 2$ and recall that the Ramsey number $r(G,H)$ is the minimum 
$n$ such that every red/blue edge-coloring of the complete $k$-uniform 
hypergraph on $n$ vertices yields either a red copy of $G$ or a blue 
copy of $H$. A 3-cycle $C_3$ is the $k$-uniform hypergraph comprising 
three edges $A, B, C$ such that every pair of them has intersection size 
1 and no point lies in all three edges. Classical results of 
Ajtai-Koml\'os-Szemer\'edi and a construction by Kim show that for 
$k=2$, we have $r(C_3, K_t)=\Theta(t^2/\log t)$, where $K_t$ is the 
complete graph on $t$ vertices. Kostochka, Mubayi, and Verstra\"ete 
proved that for $k=3$, there are positive constants $a, b$ such that
$$at^{3/2}/(\log t)^{3/4}< r(C_3, K_t) < bt^{3/2}.$$

\begin{conjecture} (Kostochka-Mubayi-Verstra\"ete)
For $k=3$, we have  $r(C_3, K_t)=o(t^{3/2})$.
\end{conjecture}

\subsection{J\o{}rgen Bang-Jensen} 

\subsubsection{Longest $(x,y)$-path in a tournament}
A digraph on at least $k+1$ vertices is $k$-strong if it remains strongly connected after the deletion of any subset $X$ of at most $k-1$ vertices. An $(x,y)$-path is a directed path from $x$ to $y$.
A digraph is {\bf hamiltonian-connected} if it contains a hamiltonian $(x,y)$-path for every choice of distinct vertices $x,y$.

\begin{theorem}\cite{thomassenJCT28}
Every 4-strong tournament is hamiltonian-connected and this is best possible.
\end{theorem} 

\begin{theorem}\cite{bangJA13}
\label{HCalg}
There exists a polynomial algorithm for deciding whether a given tournament $T$ with specified vertices $x,y$ has an $(x,y)$-hamiltonian path.
\end{theorem}

\begin{problem}[Conjecture 9.1]\cite{bangJGT28}
What is the complexity of finding the longest $(x,y)$-path in a tournament?
\end{problem}

The algorithm of Theorem \ref{HCalg} uses a divide and conquer approach to reduce a  given instance into a number of smaller instances which can either be recursively solved or for which we have a theoretical result solving the problem. Thus the approach cannot be used to solve the case where we are not looking for hamiltonian paths. 

\subsubsection{Hamiltonian paths in path-mergeable digraphs}

A digraph $D$ is {\bf path-mergeable} if, for every choice of distinct vertices $x,y\in V(D)$ and internally disjoint (only end vertices in common) $(x,y)$-paths $P_1,P_2$ there is an $(x,y)$-path $P$ in $D$ such that $V(P)=V(P_1)\cup V(P_2)$. It was shown in \cite{bangJGT20} that one can recognize path-mergeable digraphs in polynomial time. A {\bf cutvertex} in a digraph is a vertex whose removal
results in a digraph whose underlying undirected graph is disconnected. 

\begin{theorem}\cite{bangJGT20}
\label{pmdHC}
A path-mergeable digraph has a hamiltonian cycle if and only if it is strongly connected and has no cutvertex. Furthermore, a hamiltonian cycle of each block of $D$ can be produced in polynomial time.
\end{theorem}

\begin{problem}\cite{bangJGT28}
What is the complexity of the hamiltonian path problem for path-mergeable digraphs?
\end{problem}

Note that the problem is easy if $D$ is not connected or has no cutvertex so the problem is easy when the block graph of $D$ is not a path (if there is just one block the digraph has a hamiltonian cycle, by Theorem \ref{pmdHC} and if the block graph is not a path there can be no hamiltonian path. However, when the block graph is a path, the fact that we have a hamiltonian cycle in each block does not help much. In fact for every internal block with connection to its sourrounding blocks through the vertices $x,y$, we need to check the existence of an $(x,y)$-hamiltonian path.

\subsection{Klas Markstr\"om} 

Let $A$ be an $n \times n$ array with entries from $\{0,\ldots,n\}$, such that each non-zero $x$ appears in at most $n-2$ positions in $A$. ( Each entry of $A$ is just a single number.)
\begin{conjecture}
	For any $A$  there exists a latin square $L$,  using the symbols $\{1,\ldots, n\}$ such $L_{i,j}\neq A_{i,j}$
\end{conjecture}
In \cite{KO} it was proven that the conjecture holds when only two symbols appear in $A$, and a full characterization of  unavoidable arrays with two symbols was given, and a complete list of small unavoidable arrays  where each entry in $A$ can now be a list of numbers,.  In \cite{C} the conjecture was shown to hold if $n-2$ is replaced by $\frac{n}{5}$, and in \cite{O} it  was shown to hold if $A$ is a partial latin square.

\section{Problem session April 29th, 2014}

\subsection{J\o{}rgen Bang-Jensen}
A digraph $D=(V,A)$ is $k$-arc-strong if $D-A'$ remains strongly connected for every subset $A'\subseteq A$ with $|A'|\leq k-1$. We denote by $\lambda{}(D)$ the maximum $k$ such that $D$ is $k$-arc-strong.

\begin{theorem}\cite{bangJGT46} Every $k$ arc-strong tournament $T$ on $n$ vertices contains a spanning 
$k$-arc-strong subdigraph with at most $nk + 136k^2$ arcs.
\end{theorem}

\begin{theorem}\cite{bangJGT46}
Every $k$ arc-strong tournament $T=(V,A)$ on $n$ vertices contains a spanning  subdigraph $D'=(V,A')$ such that every vertex in $D'$ has in- and out-degree at least $k$ and $|A'|\leq nk+\frac{k(k-1)}{2}$ and this is best possible.
\end{theorem}

For a given tournament $T$ let $\alpha_k(T)$ denote the minimum number of arcs in a spanning subdigraph of $T$ which has minimum in- and out-degree at least $k$. For given $k$-arc-strong tournament $T$ let $\beta_k(T)$ denote the 
minimum number of arcs in a spanning $k$-arc-strong subdigraph of $T$.

\begin{conjecture}\cite{bangJGT46}
For every $k$ arc-stroing tournament $T$ we have $\alpha_k(T)=\beta_k(T)$, in particular we have 
$\beta_k(T)\leq nk+\frac{k(k-1)}{2}$.
\end{conjecture}

\begin{conjecture}\cite{bangJGT46}
There exists a polynomial algorithm for finding, in a given $k$-arc-strong tournament $T=(V,A)$ a minimum set of arcs $A'$ (of size $\beta_k(T)$) such that the subdigraph induced by $A'$ is already $k$-arc-strong.
\end{conjecture}

Note that the following theorem shows that a similar property as that conjectured above holds when we consider the minimum number $r^{arc-strong}_k(T)$ of arcs whose reversal results in a $k$-arc-strong tournament. It shows that, except when the degrees are almost right already so that some cut needs more arcs reversed, we have equality between the numbers $r^{arc-strong}_k(T)$ and $r^{deg}_k(T)$, where the later is the minimum number of arcs whose reversal in $T$ results in a tournament $T''$ with minimum in- and out-degree at least $k$.

\begin{theorem}\cite{bangDAM136}
For every tournament $T$ on at least $2k+1$ vertices the number $r^{arc-strong}_k(T)$ is equal to the maximum of the numbers
$k-\lambda{}(T)$ and $r^{deg}_k(T)$. In particular, we always have $r^{arc-strong}_k(T)\leq\frac{k(k+1)}{2}$ (equality for transitive tournaments on at least $2k+1$ vertices).
\end{theorem}

\subsection{Matas \v{S}ileikis}

A family $\mathcal{F}$ of subsets of $\left[n\right] = \left\{ 1, \dots, n \right\}$ is called
\begin{enumerate}
 \item \emph{$k$-intersecting} if for all $A,B \in \mathcal{F}$ we have $\left|A\cap B\right|\geq k$,
 \item an \emph{antichain} if for all $A,B \in \mathcal{F}$ such that $A \neq B$ we have $A \nsubseteq B$,
\end{enumerate}

 In 1964 Katona \cite{katona} (see also \cite[p. 98]{bollobas}) determined the least upper bound for the size of a $k$-intersecting family:
\begin{equation}\label{Kat}
|\mathcal{F}| \le 
\begin{cases}
  \sum_{j=t}^{n} \binom{n}{j}, \qquad  &\text{\rm if}\;\; k+n = 2t, \\
  \sum_{j=t}^{n} \binom{n}{j}+\binom{n-1}{t-1}, \qquad &\text{\rm if}\;\; k+n = 2t - 1,
\end{cases}
\end{equation}
with equality attained by the family consisting of all sets of size at least $t$ plus, when $k + n$ is odd, subsets of $[n-1]$ of size $t-1$.

In 1966 Kleitman \cite{kleitman} (see also \cite[p.~102]{bollobas}) observed that the bound \eqref{Kat} remains true under a weaker condition that $\mathcal{F}$ has diameter at most $n-k$, that is, when for every $A, B \in \mathcal{F}$ we have $|A \bigtriangleup B| \le n - k$.

In 1968 Milner \cite{milner} determined the least upper bound for the size of a $k$-intersecting antichain (which generalizes Sperner's Lemma, when $k = 0$):
  \begin{equation}\label{Mil}
    |{\cal F}| \leq 
	\binom{n}{t}, \qquad t = \left\lceil \frac {n+k} 2 \right\rceil.
  \end{equation}
\paragraph{Question.} Does the bound \eqref{Mil} still hold for antichains satisfying the weaker condition that the diameter of $\mathcal{F}$ is at most $n - k$?
\bibliographystyle{plain}

\section{Problem session May 7th, 2014}

\subsection{Imre Leader}
A Ramsey Question in the Symmetric Group.
\begin{question}
Given $k$ and $r$, does there exist $n$ such that whenever the
symmetric group $S_n$ is $k$-coloured there is a monochromatic copy of $S_r$?
\end{question}
To make sense of this, it is necessary to explain what a `copy of $S_r$'
means. We view $S_n$ as the set of all words of length $n$, on symbols
$1,...,n$, such that no symbol is repeated. Given words $x_1,...,x_r$ on
symbols $1,...,n$, such that the sum of the lengths of the $x_i$ is $n$ and no
symbol is repeated, a copy of $S_r$ means the set of all possible $r!$
concatenations (in any order) of $x_1,...,x_r$.

This is easy to check when $r=2$, but even for $r=3$ we do not know it. In fact,
we do not even know it in the case $r=3$ and $k=2$.

\end{document}